\documentclass[a4paper, 11pt]{article} 
\usepackage{amsfonts,amsmath,amssymb,amscd}
\usepackage{latexsym}
\usepackage{graphicx}
\usepackage{color} 

\newcommand{\bh}{\mathcal H}

\newcommand{\Nn}{\mathbb{N}}
\newcommand{\R}{\mathbb{R}}

\newcommand{\Zz}{\mathbb{Z}}
\newcommand{\C}{\mathbb{C}}

\newcommand{\si}{\sigma}

\newcommand{\dbar}{\bar \partial}
\newcommand{\up}{{(U,P)}}
\newcommand{\ups}{{(U_S,P_S)}}

\newcommand{\beqr}{\begin{eqnarray}}
\newcommand{\eeqr}{\end{eqnarray}}

\newcommand{\beq}{\begin{eqnarray*}}
\newcommand{\eeq}{\end{eqnarray*}}

\newcommand{\bq}{\begin{equation}}
\newcommand{\eq}{\end{equation}}

\newcommand{\bpr}{\begin{preuve}}
\newcommand{\epr}{\end{preuve}}

\newenvironment{preuve}[1][]
{\vskip 2mm  {\it \bf Proof#1. }}{$\Box$ \vskip 2mm}

\newtheorem{Theorem}{Theorem}[subsection]
\newtheorem{Remark}[Theorem]{Remark}
\newtheorem{Definition}[Theorem]{Definition}
\newtheorem{Lemma}[Theorem]{Lemma}
\newtheorem{Proposition}[Theorem]{Proposition}
\newtheorem{Corollary}[Theorem]{Corollary}

\newcommand{\equid}{\underset{d\to \infty}{\sim}}

\newcommand{\hxd}{H^0(X,L^d)}
\newcommand{\rhxd}{\R H^0(X,L^d)}

\newcommand{\rcsi}{\R C_\si}

\title{Lower estimates for the expected Betti numbers of 
\\random real hypersurfaces}
\author{Damien Gayet, Jean-Yves Welschinger}

\begin{document}
\large
\maketitle
\centerline{\textbf{Abstract}}
We estimate from below the expected Betti numbers of real hypersurfaces
taken at random in a smooth 
real projective $n$-dimensional manifold. These random hypersurfaces are chosen in the linear system of a large 
$d$-th power of a real ample line bundle equipped with a Hermitian metric of positive curvature.
As for the upper bounds that we recently established, these lower bounds
read as a product of a constant which only depends on the dimension $n$ of the manifold with
the K\"ahlerian volume of its real locus $\R X$  and  $\sqrt d^n$.
Actually, any closed affine real algebraic hypersurface appears with positive probability as part of such random 
real hypersurfaces in any ball of $\R X$ of radius $O(\frac{1}{\sqrt d}).$
\\


\section{Introduction}

What is the topology of a real hypersurface taken at random in a smooth real projective manifold?
When the latter is the projective line, this question reduces to:
how many real roots does a random real polynomial in one variable have? 
This question was answered by M. Kac \cite{Kac} in 1943 and for a different measure, by
E. Kostlan \cite{Ko}  and M. Shub and S. Smale \cite{SS} in the early 90's. In our
recent paper \cite{GaWe3} (see also \cite{GaWe2}) we did bound from above the expected Betti 
numbers of such
random real hypersurfaces in smooth real projective manifolds. Our purpose now
is to bound these Betti numbers from below, see Corollary \ref{coro 1}.

Let us first recall our framework. We denote by $X$
a smooth complex projective manifold of positive dimension $n$
defined over the reals, by $c_X: X \to X$ the induced Galois antiholomorphic involution 
 and by $\R X= \text{Fix} (c_X)$ 
the real locus of $X$ which we implicitly  assume to be non-empty.
We then consider an ample line bundle $L$ over $X$,
also defined over the reals. It comes thus equipped with an antiholomorphic 
involution $c_L: L \to L$ which turns the bundle projection map $\pi : L\to X$
into a $\Zz /2\Zz $-equivariant one, so that $c_X \circ \pi = \pi \circ c_L$.
We equip $L$ in addition with a real Hermitian metric $h$, thus 
invariant under $c_L$,
 which has a positive curvature form $\omega$ -- locally defined by $\omega = \frac{1}{2i\pi}
\partial \dbar \log h(e,e) $ for any  non-vanishing local
holomorphic section $e$ of $L$. This metric induces a K\"ahler 
metric $g_h= \omega(.\, ,i.\, )$ on $X$, which reduces to
a Riemannian metric $g_h$ on $\R X$. It then induces a $L^2$-
scalar product on every space  of global holomorphic 
real sections of tensor products
$L^d$ of $L$, $d>0$, which are 
denoted by  
$$\rhxd = \{\si \in \hxd\ | \ c_L\circ \si = \si \circ c_X\}.$$
This $L^2$-scalar product is defined by the formula
\begin{equation}\label{l2product}
 \forall (\si, \tau)\in \rhxd, \ \langle \si, \tau\rangle = \int_X h(\si,\tau)(x) dx, 
\end{equation}
where $dx$
denotes any  volume form of $X$. 
For instance, $dx$ can be chosen to be the normalized volume form $dV_h= \frac{\omega^n }{\int_X \omega^n }$.
This $L^2$-scalar 
product finally induces a Gaussian probability measure $\mu_\R$
on $\rhxd$
whose density with respect to the Lebesgue one at $\si\in \rhxd$ 
writes $\frac{1}{\sqrt\pi^{N_d}}e^{-\|\si\|^2}$,
where $N_d= \dim \hxd.$
It is with respect to this probability measure that we consider random real hypersurfaces (as
in the works \cite{Ko} and \cite{SS}) and our results hold for large values of  $d$. 

Let us recall the estimates from above that we recently established in \cite{GaWe3}.
For every $d>0$, we denote by $\R \Delta_d$ the discriminant locus 
of sections $\si \in \rhxd$ which do not vanish transversally. 
For every $\si\in \rhxd\setminus \R \Delta_d$ ,
we denote by $\rcsi = \si^{-1}(0)\cap \R X$ its smooth real vanishing locus.
Then, for every $i\in \{0, \cdots, n-1\}$, we denote by $m_i(\rcsi)$
the $i$-th Morse number of $\rcsi$ and by $b_i(\rcsi)$ its $i$-th
Betti number with real coefficients. These are defined as the infinimum over all Morse 
functions $f$ on $\rcsi$
of the number of critical points of index $i$ of $f$
and as $b_i(\rcsi)= \dim H_i(\rcsi ; \R)$
respectively. It follows from Morse theory that $b_i(\rcsi)\leq m_i(\rcsi)$
and we set 
\beq
 E(b_i) &=& \int_{\rhxd \setminus \R \Delta_d} b_i(\rcsi) d\mu_\R(\si) \\
 \eeq 
 and $
 E(m_i) = \int_{\rhxd \setminus \R \Delta_d} m_i(\rcsi) d\mu_\R(\si) .$
Then, we proved 
\begin{Theorem}[Theorem 1 of \cite{GaWe3}]\label{theo 0}
Let $X$ be a smooth real projective manifold of positive dimension $n$ and $(L,h)$ be
a real holomorphic Hermitian line bundle of positive curvature over $X$. Then, for every 
$i\in \{0, \cdots, n-1\}$, 
$$ \limsup_{d\to \infty} \frac{1}{\sqrt d^n}E(m_i) \leq \frac{1}{\sqrt \pi}e_\R(i,n-1-i) Vol_h(\R X).$$
\end{Theorem}
In Theorem \ref{theo 0},
$Vol_h(\R X)$ denotes the total Riemannian volume  of the real locus $\R X$
for the K\"ahler metric $g_h$, while $e_\R(i,n-1-i)$ is a constant 
which only depends on $i$ and the dimension of $X$. The latter originates
from random symmetric matrices and is defined as
$$e_\R(i,n-1-i)=\int_{Sym(i,n-1-i,\R )} |\det A| d\mu_\R(A),$$
where $Sym(i,n-1-i, \R) $ denotes the open cone of non-degenerated 
real symmetric matrices of size $n-1$
and signature $(i,n-1-i)$, while $d\mu_\R$ is the restriction to this cone of 
the classical Gaussian probability measure of the space of symmetric matrices,
see \cite{GaWe3}.
In particular, for fixed $i\geq 0$,
there exists $c_i>0$ such that for large values of $n$,
\begin{equation}
 \label{***}
e_\R(i,n-1-i)\leq \exp(-c_in^2),
\end{equation}
as follows from some large deviation estimates established in \cite{BAG}, see
Theorem 1.6 of \cite{GaWe3}.

Our aim now is to get similar asymptotic estimates from below
for the expected Betti numbers 
of random real hypersurfaces linearly equivalent to $L^d$,
 see Corollary \ref{coro 1}. 
These estimates will follow from our main result, Theorem \ref{theo 1},
which we now formulate.

Let $\Sigma$ be a closed hypersurface of $\R^n$,
which we do not assume to be connected. 
For every $\si \in \rhxd\setminus \R \Delta_d$,
we denote by 
$N_\Sigma(\si) $ the maximal number of disjoint open subsets 
of $\R X$ having the property that each such open subset $U'$ 
contains a hypersurface $\Sigma'$
such that $\Sigma'\subset \rcsi$
and $(U',\Sigma')$
is diffeomorphic to $(\R^n , \Sigma)$. We then set 
$$ E(N_\Sigma)= \int_{ \rhxd\setminus \R \Delta_d} N_\Sigma (\si) d\mu_\R(\si)$$
and we associate to $\Sigma$, in fact to its isotopy class in $\R^n$,
a positive constant $c_\Sigma$ out of the amount 
of transversality of a real polynomial $P$
in $n$ variables which vanishes along a hypersurface isotopic to $\Sigma$, see (\ref{csigma}).
Our main result is:
\begin{Theorem}\label{theo 1}
Let $X$ be a smooth real projective manifold of positive dimension $n$ and $(L,h)$ be
a real holomorphic Hermitian line bundle of positive curvature over $X$. Let $\Sigma$
be a closed hypersurface of $\R^n$, which does not 
need to be connected. Then, $$
\liminf_{d\to \infty} \frac{1}{\sqrt d^n} E(N_\Sigma)\geq c_\Sigma Vol_h(\R X).$$
\end{Theorem} 
In particular, when $\Sigma$ is connected, Theorem \ref{theo 1} 
bounds from below the expected number of connected 
components that are diffeomorphic to $\Sigma$
in the real vanishing locus 
of a random section $\si\in \rhxd$. As in Theorem \ref{theo 0}, the constant $c_\Sigma$ 
does not depend on the choice of
the triple $(X,L,h)$, it only depends on $\Sigma$.

Let us now denote, for every positive integer $n$, by $\bh_n$ the set of diffeormophism
classes of smooth closed  connected hypersurfaces in $\R^n$. For every $i\in \{ 0, \cdots, n-1\}$
and every $\Sigma\in \bh_n$, we denote by $b_i(\Sigma) =\dim H_i(\Sigma ;\R)$
its $i$-th Betti number with real coefficients and by $m_i(\Sigma)$
its $i$-th Morse number.
\begin{Corollary}\label{coro 1}
Let $X$ be a smooth real projective manifold of positive dimension $n$ and $(L,h)$ be
a real holomorphic Hermitian line bundle of positive curvature over $X$.
Then,
for every $i\in\{0, \cdots, n-1\}$, 
$$\liminf_{d\to \infty} \frac{1}{\sqrt d^n} E(b_i)\geq \big(\sum_{\Sigma\in \bh_n} c_\Sigma b_i(\Sigma) \big)Vol_h(\R X)$$
and likewise 
$\liminf_{d\to \infty} \frac{1}{\sqrt d^n} E(m_i)\geq \big(\sum_{\Sigma\in \bh_n} c_\Sigma m_i(\Sigma)\big) Vol_h(\R X).$
In particular, for every  $i\in\{0, \cdots, n-1\}$,
$$\liminf_{d\to \infty} \frac{1}{\sqrt d^n} E(b_i)\geq \exp\big(-2e^{70n}\big) Vol_h(\R X).$$
\end{Corollary}
The last part of Corollary \ref{coro 1} follows from the fact that for every $i\in\{0, \cdots, n-1\}$,
$\R^n$ contains the product of spheres $S^i\times S^{n-1-i}$
as a hypersurface, while $c_{S^i\times S^{n-1-i}}\geq \exp\big(-2e^{70n}\big),$
see Proposition \ref{prop 2}.
This double exponential decay has to be compared with (\ref{***})
and is  not  an optimal bound. For instance, it follows from Theorem 1 of \cite{GaWe3}
that when $n=1$, $\frac{1}{\sqrt d }E(b_0)$ converges to $ \frac{1}{\sqrt \pi } Length_h (\R X)$. 
The results given by 
Theorem \ref{theo 0} and Theorem \ref{theo 1} raise the following question:
does the quotient $\frac{E(b_i)}{Vol_h(\R X)\sqrt d^n}$ or likewise 
$\frac{E(m_i)}{Vol_h(\R X)\sqrt d^n}$ have a limit in general, which
only depends on $i\in\{0, \cdots , n-1\}$ and   the dimension $n$
of  $X$, but not on the triple $(X,L,h)$?
This holds true for $n=1$, 
see Theorem 1 of \cite{GaWe3} or also Theorem 2 of \cite{GaWe3} for similar results 
on the number of critical points of given index.

Note that another natural  Gaussian probability measure could have been chosen 
on $\rhxd$,
induced by a $L^2$-product defined by integration over $\R X$
instead of the integration over $X$ (see \S 3.1.1 of \cite{GaWe2} for a discussion on
our choice and other possible ones). This is the measure considered
by F. Nazarov and M. Sodin 
in their study of random spherical harmonics in dimension two,
see \cite{Nazarov-Sodin},
and more recently by A. Lerario and E. Lundberg in higher dimensions, see \cite{Lerario}.
The upper and lower estimates they obtain  for the number of connected components
for these spherical hamonics 
are in $d^n$ instead of $\sqrt d^n$.
These estimates are also established for homogeneous polynomials on unit
spheres in \cite{Lerario}. Note that such a behaviour was previously  guessed through computational experiments by C. Raffalli, 
while P. Sarnak and I. Wigman informed us that they were able to prove the upper estimates in $\R P^2$. 

In order to prove Theorem \ref{theo 1}, we follow the
same probability approach as Nazarov and Sodin (see \S 6.1 of \cite{Nazarov-Sodin}
or also \S 2.2 of \cite{Lerario})
which we combine with the $L^2$-estimates of H\"ormander, see
\S \ref{para II.1}. The latter 
make it possible asymptotically  to produce,
for every smooth closed hypersurface $\Sigma$ of $\R^n$
contained in a ball of radius $R$ and every 
ball $B_d$ of $\R X$ of radius $R/\sqrt d$, a section $\si \in \rhxd$ 
which vanishes transversally in $B_d$ along a hypersurface diffeomorphic to $\Sigma$,
the transversality 
being quantitative in the sense of S. Donaldson, see
Definition 7 in \cite{Donaldson} and Proposition \ref{prop 3}.
We then bound from above the expected $C^1$-norm of sections of $\rhxd$ 
in such a ball $B_d$ and deduce from Markov's inequality
that a random section in $\rhxd$ vanishes with positive probability in $B_d$ along a hypersurface 
diffeomorphic to $\Sigma$, see Proposition \ref{prop 5}. The result follows from the fact that there are more or less  $Vol_h(\R X)\sqrt d^n$ disjoint such balls in $\R X$. 
Recall that the construction in 
\cite{Gayet} of real Donaldson hypersurfaces
with many spheres in their real locus was carried out in a similar manner. 

The first part of the paper is devoted to preliminaries on closed affine real algebraic hypersurfaces
and the second one to the proofs of Theorem \ref{theo 1} and Corollary \ref{coro 1}.\\

\textit{Aknowledgements.} The research leading to these results has received funding
from the European Community's Seventh Framework Progamme 
([FP7/2007-2013] [FP7/2007-2011]) under
grant agreement $\text{n}\textsuperscript{o}$ [258204].

\section{Closed affine real algebraic hypersurfaces}

This paragraph is devoted to preliminaries. We first introduce two real functions which play a role
in the proof of Theorem \ref{theo 1}. 
Then, we associate a positive constant $c_\Sigma$ to any isotopy class of smooth closed hypersurface
$\Sigma$ of $\R^n$, see (\ref{csigma}), using a notion of regular pair given by Definition \ref{def 1}.
Finally, we estimate from below this constant in the case of product of spheres, 
see \S \ref{para exemples}.


\subsection{Two real functions}\label{para I.1}
We  introduce here two real functions $f_\tau$ and $g_R$ 
whose maximum and minimum  turn out to play a role in the proof of Theorem \ref{theo 1}.
For every $\tau>0$, we set 
\beq
 f_\tau :
 a \in  [\sqrt \tau , +\infty[ \  \mapsto \frac{1}{\sqrt \pi}(1-\frac{\tau}{a^2}) \int_a^{+\infty }e^{-t^2} dt,
 \eeq
 so that $f_\tau (\sqrt \tau) = \lim_{a\to \infty} f_\tau = 0$. 
 We set 
\bq
\label{mtau} 
  m_\tau  = \sup_{[\sqrt \tau, +\infty[} f_\tau. 
  \eq 
In particular, for every positive $\tau$, 
  \begin{eqnarray}
    m_\tau \geq f_\tau (\sqrt {\tau +1}) 
   & \geq & \frac{1}{\sqrt \pi(\tau+1)}e^{-(\sqrt{\tau+1} +1)^2}.\label{**}
 \end{eqnarray}
 The estimate (\ref{**}) is chosen in the light of the following Lemma
 \ref{lemme suivant}
 \begin{Lemma}\label{lemme suivant} For every positive $\tau$, the function $f_\tau$ 
 reaches its maximum on the interval $[\sqrt \tau, \sqrt{\tau+1}].$
 \end{Lemma}
 \bpr For every positive $\tau$ and every $a\geq \sqrt \tau $,
 $ f'_\tau (a) = \frac{1}{\sqrt \pi}\left(\frac{2\tau}{a^3}\int_a^{+\infty }e^{-t^2} dt - (1-\frac{\tau}{a^2})e^{-a^2}\right),
 $
 so that if $f_\tau$ reaches its maximum at the point $c\in [\sqrt \tau, +\infty[$,
$ \frac{2\tau }{c^3}\int_c^{+\infty }e^{-t} dt = (1-\frac{\tau}{c^2}) e^{-c^2}.$
Now $\int_c^{+\infty }e^{-t^2} dt \leq \frac{e^{-c^2}}{2c}$, so that
 $(1-\frac{\tau}{c^2})\leq \frac{\tau}{c^4}$ and $c^2-\tau \leq \frac{\tau}{c^2}\leq 1$.
 Hence $c\leq \sqrt{1+\tau}.$
 \epr
 
Likewise, for every positive $R$ and every positive integer $n$, we set
\beq
g_R :
s \in \R^*_+ \mapsto \frac{(R+s)^{2n}}{s^{2n}} e^{\pi (R+s)^2},
\eeq
so that $\lim_{s\to 0} g_R(s) = \lim_{s\to +\infty } g_R(s)= +\infty.$
We set 
\bq \label{rhoR}
\rho_R = \inf_{\R^*_+ } g_R.
\eq 
In particular, for every positive $R$, 
\bq\label{ineq rhoR}
\rho_R \leq g_R (R) = 4^n \exp(4\pi R^2).
\eq

\subsection{Real polynomials and transversality}\label{para I.2}

We introduce here the notion of regular pair, see Definition \ref{def 1},
and the constant $c_\Sigma$ associated to any isotopy class of smooth closed hypersurface $\Sigma$
of $\R^n$, see (\ref{csigma}).
%
\begin{Lemma}\label{lemme 2}
If $P = \sum_{(i_1, \cdots, i_n)\in \Nn^n} a_{i_1, \cdots, i_n} z^{i_1}_1\cdots z_n^{i_n} \in \R [X_1, \cdots, X_n]$,
then
$$ \| P\|^2_{L^2} = \int_{\C^n} |P(z)|^2 e^{-\pi \|z\|^2} dz = \sum_{(i_1, \cdots, i_n)\in \Nn^n} |a_{i_1, \cdots, i_n} |^2 \frac{i_1!\cdots i_n !}{\pi^{i_1+\cdots +i_n}}.$$
\end{Lemma}
\bpr
We note that $\|P\|_{L^2}^2 = \sum_{I,J\in \Nn^n} a_I \overline{ a_J } \int_{\C^n} z^I\bar z^J e^{-\pi\|z\|^2} dz.$
But for every $k\not=0, \ \int_\C z^k e^{-\pi|z|^2}dz = \int_\C \bar z^k e^{-\pi|z|^2}dz =0,$
whereas for every non-negative $k$,
$$ \int_\C |z|^{2k} e^{-\pi|z|^2} dz = 2\pi\int_0^{+\infty} r^{2k+1} e^{-\pi r^2} dr = \pi\int_0^{+\infty } t^k e^{-\pi t} dt = 
\frac{1}{\pi^k} \Gamma(k+1) = \frac{k!}{\pi^k}.$$
The result follows then from Fubini's Theorem.
\epr
\begin{Definition}\label{def 1}
Let $U$ be a bounded open subset of $\R^n$ and $P\in \R[X_1, \cdots X_n]$, $n>0$. 
The pair $(U,P)$ is said to be regular if and only if 
\begin{enumerate}
\item zero is a regular value of the restriction of $P$ to  $U$,
\item the vanishing locus of $P$ in $U$ is compact.
\end{enumerate}
\end{Definition}

\begin{Definition}\label{def 2}
For every regular pair $(U,P) $ given by Definition \ref{def 1}, 
we denote by $\mathcal T_{(U,P)} $ the set of $(\delta, \epsilon)\in (\R^*_+)^2$ such that
\begin{enumerate}
\item there exists a compact subset $K$ of $U$ satisfying $\inf_{U\setminus K} |P| >\delta$,
\item for every $y\in U$, $|P(y)| <\delta \Rightarrow \| d_{|y} P\| >\epsilon$,
where $\|d_{|y} P\|^2= \sum_{i=1}^n |\frac{\partial P}{\partial x_i}|^2.$
\end{enumerate}
\end{Definition}
We then set for every regular pair $\up$, $R_{(U,P)} = \max (1, \sup_{y\in U} \|y\|) $ and 
\bq \label{tauUP}
\tau_{(U,P)} = 2\rho_{R_{(U,P)}} \|P\|^2_{L^2} \inf_{(\delta,\epsilon)\in {\mathcal T_{(U,P)}}} (\frac{1}{\delta^2}+ \frac{\pi n}{\epsilon^2})\in \R^*_+,
\eq
where $\rho_{R_\up}$ is defined by (\ref{rhoR}).

Now, let $\Sigma$ be a closed hypersurface  of $\R^n$, not necessarily connected. 
We denote by $\mathcal I _\Sigma$ the set of regular pairs $(U,P)$ given by Definition \ref{def 1}, such that
the vanishing locus of $P$ in $U$ contains a subset isotopic to $\Sigma$ in $\R^n$.
It follows from Nash's Theorem for hypersurfaces in $\R^n$ that 
$\mathcal I_\Sigma$ is non empty, see \cite{Nash}. 
We  then set
\bq\label{csigma}
c_\Sigma = \sup_{(U,P)\in \mathcal I_\Sigma} \left(\frac{m_{\tau_{(U,P)}}}{2^n Vol(B(R_{(U,P)}))}\right),
\eq
where $m_{\tau_\up}$ is defined by (\ref{mtau}) and $Vol(B(R_{(U,P)}))$ denotes the  volume of the Euclidean ball of radius $R_{(U,P)}$ in $\R^n$. 
From (\ref{**}) follows that for every $(U,P)\in \mathcal I_\Sigma$,
\beqr
c_\Sigma 
& \geq & \frac{\lfloor \frac{n}{2}\rfloor ! \exp\big(-(\sqrt{\tau_\up +1} +1)^2\big)}{2^{n+1}\pi^{\lfloor \frac{n}{2}\rfloor } R^n_\up (1+\tau_\up)\sqrt \pi},\label{ineq csigma}
\eeqr
since the volume of the ball of radius $R_\up$ in $\R^n$ is bounded from below by $\frac{2 \pi^{\lfloor \frac{n}{2} \rfloor} }{\lfloor \frac{n}{2} \rfloor!}R^n_\up$
for every $n>0$.
For large values of $\tau_\up$, as the ones which appear in the examples given in \S \ref{para exemples},
we deduce from (\ref{ineq csigma}) that 
\begin{equation}\label{8bis}
c_\Sigma \geq \exp(-2\tau_\up).
\end{equation}
\subsection{Examples}\label{para exemples}
\subsubsection{The spheres}
\begin{Proposition}\label{prop 1}
For every positive integer  $n$, $c_{S^{n-1}}\geq \exp (-2e^{43n}).$
\end{Proposition}
For every $n>0$, we set $P_S(x_1, \cdots, x_n) = \sum_{j=1}^n x_j^2 - \sqrt n - 1$
and  $ U_S = \{(x_1, \cdots , x_n)\in \R^n \, | \, \sum_{j=1}^n x_j^2 < \sqrt n +2\}.$
The pair $\ups$ is regular in the sense of Definition \ref{def 1} and $P^{-1}_S(0) \subset U_S$
is isotopic in $\R^n$ to the unit sphere $S^{n-1}$. 
\begin{Lemma} \label{lemme 3}
For every $n>0$ and every $0<\delta<1$, $\big(\delta, 2\sqrt{\sqrt n +1 -\delta}\big)\in \mathcal T_{\ups}.$
\end{Lemma}
\bpr For every $x\in \R^n$ and $\delta>0$,
\beq 
|P_S(x)| <\delta &\Leftrightarrow & \sqrt n +1 -\delta < \|x\|^2< \sqrt n +1+\delta \\
& \Rightarrow & \|d_{|x} P_S\|^2 = 4\|x\|^2 > 4(\sqrt n +1-\delta).
\eeq
Moreover, when $0<\delta <1$, $K_\delta = \{x\in U \ | \ \sqrt n +1 - \delta\leq \|x\|^2\leq \sqrt n +1+\delta\}$
is compact in $U_S$. We deduce  that $(\delta,\epsilon)\in \mathcal T_\ups$ for
$\epsilon^2 = 4(\sqrt n +1-\delta)$.
\epr
\bpr [ of Proposition \ref{prop 1}]
For every positive integer  $n$, $R^2_\ups = \sqrt n +2$, while from Lemma \ref{lemme 2}, 
$\|P_S\|^2_{L^2} = (\sqrt n +1)^2 + \frac{2n}{\pi^2} \leq 5 n.$
From (\ref{tauUP}) and Lemma \ref{lemme 3},
we deduce
\beq
\tau_\ups & \leq & \rho_{R_{(U_S,P_S)}} 10 n (1+\frac{\pi \sqrt n }{4})\\
& \leq & \exp\left(n\ln 4 + 4\pi(\sqrt n +2) + \ln \big(10n + \frac{\pi \sqrt n }{4}\big)\right) \text{ by (\ref{ineq rhoR})}\\
& \leq & \exp(43 n).
\eeq
The estimate
$c_{S^{n-1}}  \geq \exp (-2e^{43n}) $ follows then from (\ref{8bis}).
\epr
\subsubsection{Products of spheres}
\begin{Proposition}\label{prop 2}
 For every positive integer $n$ and every $0\leq i \leq n-1$,
$c_{S^i \times S^{n-i-1} }\geq \exp \big(-2 e^{70 n} \big).$
\end{Proposition}
For every $n>0$ and every $0\leq i\leq n-1$, we set 
$$
 Q_i 
 ((x_1, \cdots, x_{i+1}), (y_1, \cdots, y_{n-i-1}))  =  \big(\sum_{j=1}^{i+1} x_j^2 -2\big)^2 + \sum_{j=1}^{n-i-1} y^2_j - 1
$$
and $ U_{Q_i} = \{(x,y)\in \R^{i+1} \times \R^{n-1-i} \, | \, \|x\|^2 + \|y\|^2 <5\}.$
The pair $(U_{Q_i}, Q_i)$ is regular in the sense of Definition \ref{def 1} and $Q_i^{-1}(0) \subset V_{Q_i}$
is isotopic in $\R^n$ to the product $S^i \times S^{n-i-1}$ of the unit spheres in $\R^{i+1} $ and $\R^{n-i}$.

\begin{Lemma} \label{lemme 4}
For every positive integer $n$ and  every $0\leq i\leq n-1$, $$\big(\frac{1}{2\sqrt n}, 2\sqrt{1-\frac{1}{2\sqrt n}}\big) \in \mathcal T_{(U_{Q_i}, Q_i)}.$$
\end{Lemma}
\bpr For every $(x,y)\in \R^{i+1}\times \R^{n-i-1}$ and every $0<\delta <\frac{1}{2}$,
\beq
|Q_i(x,y)| <\delta &\Leftrightarrow & 1-\delta < (\|x\|^2-2)^2 + \|y\|^2<1+\delta\\
& \Rightarrow & \|d_{|(x,y)} Q_i\|^2 = 4\|y\|^2 + 16 \|x\|^2 (\|x\|^2-2)^2
\eeq
with $\|x\|^2 >2 - \sqrt{1+\delta} > 1/2$. Thus, $\|d_{|(x,y)} Q_i\|^2 > 4 \big((\|x\|^2-2)^2 + \|y\|^2\big) >4(1-\delta)$ and
we deduce the result by choosing $\delta = \frac{1}{2\sqrt n}$.
\epr
\bpr[ of Proposition \ref{prop 2}] For every positive integer $n$ and every $0\leq i \leq n-1$, 
$R^2_{(U_{Q_i}, Q_i)} = 5$,
$\rho_{R_{(U_{Q_i}, Q_i)}}\leq 4^n \exp(20\pi) $ by (\ref{ineq rhoR}), while from Lemma \ref{lemme 2},
$$ \| Q_i\|^2_{L^2} = 9 + \frac{2}{\pi^2}(n-i-1)+\frac{32}{\pi^2}(i+1) +\frac{24}{\pi^4}(i+1) + \frac{16}{\pi^4} {i+1\choose 2}\leq 13 n^2.$$
We deduce from (\ref{tauUP}) and  Lemma \ref{lemme 4} the  upper estimate
$\tau_{(U_{Q_i}, Q_i)}  \leq  156 n^3 4^n e^{20\pi} 
 \leq  e^{70n} $
since $\ln n \leq n-1$.
We then deduce from (\ref{8bis}) the lower estimate
$c_{S^i \times S^{n-i-1}}  \geq  \exp(-2e^{70n}).$
\epr
\begin{Remark} The lower estimates given by Propositions \ref{prop 1} and  \ref{prop 2}
are far from being optimal.
\end{Remark} 

\section{Lower bounds for the Betti numbers of random real algebraic hypersurfaces}

We first implement the affine real algebraic hypersurfaces in every smooth real projective manifold at the scale $1/\sqrt d$ 
 thanks to H\"ormander $L^2$-estimates, see Proposition \ref{prop 3}. We then follow the approach of
F. Nazarov and M. Sodin (see \cite{Nazarov-Sodin} or also \cite{Lerario})
by first estimating the expected local $C^1$-norm of sections, 
see Proposition \ref{prop 4}, and then deducing a positive probability of presence of such
affine real algebraic hypersurfaces in the vanishing locus of random sections
in any ball of radius $O(\frac{1}{\sqrt d})$, see Proposition 
\ref{prop 5}. Theorem \ref{theo 1} and Corollary \ref{coro 1} follow.

\subsection{H\"ormander sections}\label{para II.1}

\begin{Definition}\label{trivi}
Let $(X,c_X)$ be a smooth real projective manifold of positive dimension $n$ and $(L,h,c_L)$ be
a real holomorphic Hermitian line bundle of positive curvature over $X$. For every $x$ in $ \R X$, let us
call a $h$-$trivialization$ of $L$ in the neighbourhood of $x$ the following data:
\begin{enumerate}
\item a local holomorphic chart $\psi_x : (W_x,x) \subset X \to (V_x,0)\subset \C^n$
such that \label{condition 1}
\begin{enumerate}
\item $\psi_x\circ c_X = conj \circ \psi_x$, where $conj:  (y_1, \cdots, y_n)\in \C^n \mapsto  (\bar y_1, \cdots, \bar y_n) \in \C^n$
\item $d_{|x}\psi_x : (T_x X, g_h)\to \C^n $ be an isometry
\end{enumerate}
\item a non-vanishing holomorphic section $e$ of $L$ defined over  $W_x$ and such that \label{condition 2}
\begin{enumerate}
\item  $c_L\circ e\circ c_X=e$
\item $\phi = -\log h(e,e)$ vanishes at $x$ and is positive everywhere else
\item there exist positive constants $\alpha_1$, $\alpha_2$ such that on $V_x$.
\begin{equation}\label{phi}
 \alpha_1 \|y\|^3 \leq \phi\circ \psi_x^{-1}(y)- \pi\|y\|^2 \leq \alpha_2 \| y\|^3,
 \end{equation}
\end{enumerate}
\end{enumerate}
\end{Definition}

\begin{Definition}\label{def htrivi}
Let $X$ be a smooth real projective manifold of positive dimension $n$ and $(L,h)$ be
a real holomorphic Hermitian line bundle of positive curvature over $X$.
A  {\it field of $h$-trivializations} on $\R X$ is the data at every real point $x$ of a $h$-trivialization
in the neighbourhood of $x$ such that 
the open subset $V_x$ of $\C^n$  given
by Definition \ref{trivi} does not depend on $x\in \R X$
and such that the composition $\phi\circ \psi_x^{-1}$ gets uniformly bounded from below 
by a positive constant on this open set $V=V_x$, while the constants $\alpha_1$, $\alpha_2$  
can be chosen not to depend on $x\in \R X$.
\end{Definition}

\begin{Lemma}
Any smooth real projective manifold of positive dimension equipped with a real holomorphic
Hermitian bundle of positive curvature admits a field of $h$-trivializations.
\end{Lemma}
\bpr
Let $(X,c_X)$ be a smooth real projective manifold of positive dimension $n$ and $(L,h,c_L)$ be
a real holomorphic Hermitian line bundle  with positive curvature $\omega$ over $X$.
Let $x$ be a real point of $X$ and let us first prove the existence of a $h$-trivialization
of $L$ near $x$. The existence of the local chart $\psi_x : (W_x,x)\to (V_x,0)$
satisfying the first condition  of Definition \ref{trivi} is given 
by definition. Now, restricting $W_x$ if necessary and averaging a local holomorphic section,
we get a local holomorphic section $\tilde e$ of $L_{|W_x}$ which does not vanish and
satisfies $c_L \circ \tilde e = \tilde e \circ c_X$. The plurisubharmonic function
$\tilde \phi =- \log h(\tilde e,\tilde e)$  
takes real values. Its composition $\tilde \phi \circ \psi_x^{-1}$ writes
$ \tilde \phi \circ \psi_x^{-1}= \Re  \phi_1 +  \phi_2$, where $\phi_1\in \C[y_1, \cdots, y_n]$ 
is a degree  two polynomial and $\phi_2(y) = \pi \|y\| ^2 + O(\|y\|^3)$, since
the Hermitian part of the second derivative of $\tilde \phi$ at $x$ is $\pi g_h$ by definition.
We then set, following \cite{Hormander1}, $e = \exp ( \phi_1 \circ \psi_x)\tilde e $ which satisfies the second condition
of Definition \ref{trivi} after restricting the open subset $W_x$ if necessary. 

How small has to be chosen $W_x$ depends on the higher order derivatives of $ \phi_2$.
However, these higher order derivatives are the same as the ones 
of $\tilde \phi \circ \psi_x^{-1}$ since
they are not affected by $\phi_1$. Now, we can cover $\R X$ with the supports 
of finitely many real sections $\tilde e_1, \cdots , \tilde e_k$. The derivative of these sections
are uniformly bounded over $\R X$. We can thus choose a $h$-trivialization
near every point $x$ of $\R X$ in such a way that the open subset
$V_x$ of $\C^n$ does not depend on $x\in \R X $. Restricting $V_x$
if necessary, this ensures the existence of a field of $h$-trivialization on $\R X$.
\epr

For every positive $d$ and every  $\sigma \in \rhxd$, 
  $\| \si\|_{L^2(h)}$ denotes the $L^2$-norm
  for the normalized volume form induced by the K\"ahler form $\omega$,
  that is 
  $$ \| \si\|_{L^2(h)}^2  = \int_X \|\si\|^2_{h^d} dV_h,$$ 
  where $dV_h= \frac{\omega^n  }{\int_X \omega^n }.$
 Moreover, if
the restriction of $\sigma$ to $W_x$
writes $\sigma = f_\sigma e^d$ for some holomorphic  $f_\sigma : W_x\to \C $, 
we  set $|\sigma| = |f_\si|$ 
so that on $W_x$,
$ \| \sigma \|^2_{h_d} = |\sigma|^2 \exp(-d\phi )$,
and for every $z$ in $W_x$,
\begin{equation}\label{derivee}
|d_{|z}\sigma| = |d_{|y}(f_\si\circ \psi_x^{-1})|, \text{ where } y=\psi_x(z).
\end{equation}
We also denote, for every small enough $R>0$, by $B(x,R)\subset W_x$
the ball centered at $x$ and of radius $R$ for the flat metric of $V$
pulled back by $\psi_x$, so that  $B(x,R)= \psi_x^{-1}(B(0,R))$. We finally denote by $\delta_L = \int_X c_1(L)^n $ 
  the degree of the bundle $L$ and recall that the notion of regular pair has been defined
in Definition \ref{def 1}. 
\begin{Proposition}\label{prop 3}
Let $X$ be a smooth real projective manifold of positive dimension $n$ and $(L,h)$ be
a real holomorphic Hermitian line bundle  with positive curvature over $X$.
We choose a field of $h$-trivializations on $\R X$.  Then, for every regular pair $(U,P)$, every  large enough integer  $d$ 
 and every  $x$ in $\R X$, there exist $\sigma_\up\in \rhxd $   and an open subset $U_d\subset B(x,\frac{R_\up}{\sqrt d})$
such that 
\begin{enumerate}
\item $ \|\si_\up\|_{L^2(h)}$ be equivalent to $\frac{\|P\|_{L^2} }{\sqrt{\delta_L}}$ as $d$ grows to infinity \label{un}
\item $ (U_d, \sigma_\up^{-1} (0) \cap U_d) $ be diffeomorphic  to $(U, P^{-1}(0)\cap U) \subset \R^n$
\item for every $(\delta,\epsilon)\in \mathcal T_\up$ given by Definition \ref{def 2},  there exists a compact subset $K_d\subset U_d$
such that $$\inf_{U_d\setminus K_d} |\si_\up| >\frac{\delta}{2}\sqrt d^n, $$
while for every $y$ in $U_d$,  
\begin{equation}\label{donaldson}
 |\si_\up (y)|< \frac{\delta}{2} \sqrt d^n \Rightarrow |d_{|y}\si_\up|>\frac{\epsilon}{2}\sqrt d^{n+1}.
\end{equation}
\end{enumerate}
\end{Proposition}
Recall that  the norm of the derivative is given by (\ref{derivee}), and note that the quantitative transversality condition  (\ref{donaldson}) is the 
one used by Donaldson in \cite{Donaldson}.\\

Under the hypotheses of Proposition \ref{prop 3}, let $x\in \R X$. We set, for every  large enough $d$, 
$U_d = \psi_x^{-1}(\frac{1}{\sqrt d} U) \subset B(x,\frac{R_\up}{\sqrt d}).$
Let $\chi : \C^n \to [0,1]$ be a smooth function with compact support
in $V$ which equals one in a neighbourhood of the origin. 
 Then, let $\sigma$ be the global smooth section of $L^d$ defined by 
$\sigma_{|(X\setminus W_x)} = 0$ and 
$$ \ \sigma_{|W_x}= (\chi \circ \psi_x) P(\sqrt d \psi_x) e^d.$$
From  the $L^2$-estimates of H\"ormander, see  \cite{Hormander1}, \cite{Skoda}, \cite{Tian},  there exists a global section 
$\tau$ of $L^d$  such that $\dbar \tau = \dbar \sigma $ and $\|\tau\|_{L^2(h)} \leq \|\dbar \sigma\|_{L^2(h)}$ for $d$ 
large enough.
This section can be chosen orthogonal to the holomorphic sections and is then unique, in particular real, 
so that $c_{L^d} \circ \tau \circ c_X = \tau. $
Moreover,
\begin{Lemma}\label{estimees}
There exist positive constants $c_1$ and $c_2$ which do not depend on  $x\in \R X$ and satisfy 
$\|\tau\|_{L^2(h)}\leq c_1e^{-c_2d} $ as well as 
$|\tau|_{C^1(U_d)} \leq c_1e^{-c_2d},$  where the $C^1$-norm is defined by (\ref{derivee}).
\end{Lemma}
\bpr The $L^2$-estimates of H\"ormander read for large enough $d$ 
\beqr\label{inequ}
\|\tau\|^2_{L^2(h)}&\leq & \int_X \|\dbar \sigma\|^2 dV_h
 =  \int_{W_x} | \psi_x^*\dbar \chi |^2 |P(  \sqrt d \psi_x )|^2 e^{-d\phi} dV_h\\
& \leq &  d^{deg(P)} \sup_{V}\left( |\dbar \chi|^2 \left|\frac{P(\sqrt d \, .\,)}{\sqrt d ^{deg P}}\right|^2\right) 
e^{-d \underset{supp \, (\dbar \chi)}{\inf}( \phi \circ \psi_x^{-1})},
\eeqr
so that there exist positive constants $c_1, \tilde c_2$, not depending on $x$, such that $\|\tau\|_{L^2(h)} \leq c_1 \exp(- \tilde c_2d)$.
Now, since $\tau \circ \psi_x^{-1}$ is holomorphic on $\chi^{-1}(1)$, the mean value 
inequality for plurisubharmonic functions implies that for every $z$ in $U_d$,
\beq
 |\tau(z)|^2 &\leq &\frac{1}{Vol(B(\frac{1}{\sqrt d}))}\int_{B(\psi_x(z),\frac{1}{\sqrt d}))}|\tau\circ \psi_x^{-1}(y)|^2 dy \\
 & \leq &  \frac{1  }{Vol(B(\frac{1}{\sqrt d}))}
 \int_{B(z,\frac{1}{\sqrt d}))}\|\tau\|_{h^d}^2e^{d\phi}  \psi_x^*dy\\
 & \leq &  \frac{1 }{Vol(B(\frac{1}{\sqrt d}))}\sup_{B(z,\frac{1}{\sqrt d}))}  \big( e^{d\phi  } |\det d_{|z}\psi_x| \big)\|\tau\|_{L^2(h)}^2,
 \eeq 
 where
 the determinant $|\det d_{|z}\psi_x|$ is computed with respect to the volume forms
 $dV_h$ and $dy$. We deduce that there exists a constant
$c_3$ not depending on  $x\in \R X$ such that
 $|\tau|^2\leq c_3 d^n e^{-2 \tilde c_2 d}$ on $U_d$,
  see (\ref {phi}) and Definition \ref{def htrivi}.
The estimate for $|d\tau|$ is proved along the same lines. 
\epr
\bpr[ of Proposition \ref{prop 3}]
We set $\si_\up = \sqrt d^n (\si- \tau) $ and $K_d = \psi_x^{-1}(\frac{1}{\sqrt d}K),$ 
see Definition \ref{def 2}. The section $\si_\up$ is global and holomorphic.
Lemma \ref{estimees} shows that on $U_d$,  $\si_\up$ is
a small perturbation of $\sqrt d^n \si$.
In particular,
\beqr
\| \sigma_\up\|_{L^2(h)}^2 
\equid  d^n \int_{\chi^{-1} (1)} |P(\sqrt d y)|^2 e^{-d \pi \|y\|^2} \frac{dy}{\delta_L}
 \underset{d\to \infty}{\to}  \frac{\| P\|^2_{L^2}}{\delta_L}.
\eeqr
Moreover, for every pair $(\delta, \epsilon)\in \mathcal T_\up$ and every $z\in U_d\setminus K_d$,
\beqr
|\frac{1}{\sqrt d^n} \si_\up (z)| &=& |\sigma(z) - \tau (z)|\\
& \geq & | \sigma(z)| - \sup_{B\big(x,\frac{R_\up}{\sqrt d}\big)} |\tau| \\
& \geq & |P( \sqrt d \psi_x( z))| -  \sup_{B\big(x,\frac{R_\up}{\sqrt d}\big)} |\tau| \\
&>& \delta -  \sup_{B\big(x,\frac{R_\up}{\sqrt d}\big)} |\tau| \text{ from Definition \ref{def 2}.}
\eeqr

 Thus, by Lemma \ref{estimees}, if $d$ is large enough, $\inf_{U_d\setminus K_d} |\sigma_\up| >\frac{\delta}{2}\sqrt d^n$
 whenever $x\in \R X$.
Moreover, for every $z\in U_d$, 
\beq
|\sigma_\up(z)| < \frac{\delta}{2}\sqrt d^n &\Rightarrow &|\sigma(z) - \tau(z)| <\frac{\delta}{2}\\
& \Rightarrow & |\sigma(z)| < \frac{\delta}{2} + |\tau(z)| \\
& \Rightarrow & |P(\sqrt d \psi_x (z)| < \delta 
\eeq
for $d$ large enough, whatever  $x\in \R X$ is.
Thus,
\beq
|\sigma_\up(z)| < \frac{\delta}{2}\sqrt d^n &\Rightarrow & |d_{|\sqrt d \psi_x (z)}P | >\epsilon \\
& \Rightarrow & |d_{|z}\si | >\epsilon \sqrt d \text{ using  notation (\ref{derivee})}\\
& \Rightarrow & | d_{|z} \sigma_\up| > \frac{\epsilon}{2} \sqrt d^{n+1} 
\eeq
for $d$ large enough by Lemma \ref{estimees}. 
Finally, Lemma \ref{estimees} together with 
 Lemma \ref{lemme 5} imply that  $(\si_\up^{-1}(0)\cap U_d)$ is isotopic to $(\si^{-1}(0)\cap U_d)$
and so diffeomorphic  to $(P^{-1}(0)\cap U)$  when $d$ is large enough.
\epr
\begin{Lemma}\label{lemme 5}
Let $U$ be an open subset of $\R^n $, $f: U\to \R$ be a  function of class $C^1$  and $(\delta, \epsilon)\in (\R^*_+)^2$
be such that 
\begin{enumerate}
\item there exists a compact subset $K$ of $ U$ such that $\inf_{U\setminus K} |f| >\delta$,
\item for every $ y$ in $U$, $|f(y)|< \delta \Rightarrow |df_{|y} |>\epsilon.$
\end{enumerate}
Then, for every function $g : U\to \R$ of class $C^1$ such that $\sup_U|g| <\delta$ and $\sup_U |dg| <\epsilon$,
zero is a regular value of $f+g$ and $(f+g)^{-1}(0)$ is compact and isotopic to $f^{-1}(0)$ in $U$.
\end{Lemma}
\bpr
For every $t\in [0,1]$ and every $y\in U$, 
$$f+tg(y) = 0 \Rightarrow |f(y)| = |tg(y)| < \delta.$$
The point $y$ is then contained in $K$ and $|df_{|y}|>\epsilon$.
Hence, we have
$$ |d_{|y} (f+tg)| \geq |d_{|y} f| - |td_{|y} g|  >0,$$
so that $0$ is a regular value of $f+tg$. The hypersurface $\Sigma_t = (f+tg)^{-1}(0)$
is smooth and included in $K$ from the implicit function Theorem. It
produces an isotopy between $f^{-1}(0) $ and $(f+g)^{-1}(0)$. 
\epr
\subsection{The expected local $C^1$-norm of sections}
The following Proposition \ref{prop 4}
computes the expected local $C^1$-norm of sections. It is inspired by an analogous result  of 
F. Nazarov and M. Sodin, see \cite{Nazarov-Sodin} (or also \cite{Lerario}).
Recall that we denote by $\delta_L$ the degree $\int_X c_1(L)^n $ of the line bundle
$L$ over $X$,  that $| . | $ denotes the modulus evaluated in the charts given
by $h$-trivializations, see (\ref{derivee}), and  that the constant $\rho_R$ is defined by (\ref{rhoR}).
Finally, we denote by  $v$ the density of $dV_h$ with respect to the 
 volume form $dx$ chosen in  (\ref{l2product}) to define the 
$L^2$-product, so that 
$dV_h = v(x)dx$.

\begin{Proposition}\label{prop 4}
Let $X$ be a smooth real projective manifold of positive dimension $n$ and $(L,h)$ be
a real holomorphic Hermitian line bundle of positive curvature over $X$. We equip 
$\R X$ with a field of $h$-trivializations. Then, for every positive  $R$,
\beq
\limsup_{d\to \infty} \sup_{x\in \R X}\frac{1}{d^n} E\big(\sup_{B(x,\frac{R}{\sqrt d})} \frac{|\si|^2}{v(x)}\big) 
& \leq & \frac{1}{2} \delta_L  \rho_R  \text{ and}\\
\limsup_{d\to \infty} \sup_{x\in \R X} \frac{1}{d^{n+1}} E\big(\sup_{B(x,\frac{R}{\sqrt d})} \frac{|d\si|^2}{v(x)}\big) 
& \leq & \frac{1}{2} \pi n \delta_L  \rho_R.
\eeq
\end{Proposition}
\bpr
Let $R>0$, $x\in \R X$ and $W_x$ be a neighbourhood given by the $h$-trivialization.
When $d$ is large enough, $B(x,R/\sqrt d)\subset W_x$
and $\phi\circ \psi_x^{-1}(y)= \pi \|y\|^2 + o(\|y\|^2).$
We deduce from the mean value inequality that for every $s\in \R_+$ and $\si \in \rhxd$,
\beq
\forall z\in B(x,\frac{R}{\sqrt d}), \, |\si(z)|^2 & \leq & \frac{1}{Vol(B(\frac{s}{\sqrt d}))}\int_{B(z,\frac{s}{\sqrt d})} |\si|^2 \psi_x^*dy\\
& \leq & \frac{1}{Vol(B(\frac{s}{\sqrt d}))}\int_{B(x,\frac{R+s}{\sqrt d})}  |\si|^2 \psi_x^*dy\\
\eeq
Thus, 
$ \sup_{B(x,\frac{R}{\sqrt d})}  |\si|^2\leq \frac{1}{Vol(B(\frac{s}{\sqrt d}))}\int_{B(x,\frac{R+s}{\sqrt d})} |\si|^2 \psi_x^*dy$
and after exchange of the integrals,
$$E \big(\sup_{B(x,\frac{R}{\sqrt d})}|\si|^2 \big)\leq \frac{1}{Vol(B(\frac{s}{\sqrt d}))}\int_{B(x,\frac{R+s}{\sqrt d})} E( |\si|^2 ) \psi_x^*dy.$$
Then, let $z\in B(x,\frac{R+s}{\sqrt d})$ and $\si_0\in \rhxd$
be the Bergman section at $z$. Its norm equals one and it is orthogonal  to the space
of sections  vanishing at $z$. 
Assume for the moment that the volume form $dx$ chosen to define the $L^2$-scalar product 
equals $dV_h$, so that $v=1$. 
Then, from Lemma 2.2 of \cite{Tian}
(see also \cite{Bouche}, \cite{GaWe2}),
$$ \|\si_0(z)\|^2_{h^d}\equid \delta_L d^n.$$
But $\|\si_0(z)\|^2_{h^d}= |\si_0(z)|^2 e^{-d\phi(z)},$
from which we deduce 
$$ |\si_0(z)|^2\leq \delta_L d^n e^{\pi(R+s)^2}+ o(d^n),$$
where the $o(d^n)$ can be chosen not to depend on $x\in \R X$.
As a consequence,  
\beq
E(|\si(z)|^2)&=& \int_{\rhxd} |\si(z)|^2 d\mu_\R (\si) =  |\si_0(z)|^2\frac{1}{\sqrt \pi} \int_\R a^2 e^{-a^2} da\\
& \leq & \frac{1}{2} \delta_L d^n e^{\pi(R+s)^2} + o(d^n),
\eeq
and finally, 
\beq
\sup_{x\in \R X}E\big(\sup_{B(x,\frac{R}{\sqrt d})} |\si|^2\big)&\leq & \frac{\delta_L d^n}{2Vol(B(\frac{s}{\sqrt d}))}\int_{B(0,\frac{R+s}{\sqrt d})} e^{\pi(R+s)^2} dy +o(d^n)\\
& \leq & \frac{1}{2} \delta_L d^n e^{\pi(R+s)^2} \frac{(R+s)^{2n}}{s^{2n}}+ o(d^n).
\eeq
Choosing $s\in \R^*_+$ such that $g_R(s)=\rho_R$, see (\ref{rhoR}),
we deduce that 
$$\limsup_{d\to \infty} \sup_{x\in \R X}\frac{1}{d^n}E(\sup_{B(x,\frac{R}{\sqrt d})} |\si|^2) \leq \frac{1}{2} \rho_R \delta_L .$$
Likewise, we deduce from the mean value inequality that  for every $s\in \R_+$,
$j\in\{1, \cdots, n\}$ and $z\in B(x,\frac{R}{\sqrt d})$,
\beq
\big|\frac{\partial \si}{\partial y_j}\big|^2(z)
& \leq &  \frac{1}{Vol\big(B(\frac{s}{\sqrt d })\big)}\int_{B(0,\frac{s+R}{\sqrt d} )}\big|\frac{\partial (\si\circ\psi_x^{-1})}{\partial y_j}\big|^2(y) dy,
\eeq
from which follows after summation over $j\in \{1, \cdots, n\}$ that
$$ E\big(\sup_{B(x,\frac{R}{\sqrt d}) }|d\si|^2\big)\leq \frac{1}{Vol\big(B(\frac{s}{\sqrt d } )\big)}\int_{B(0,\frac{s+R}{\sqrt d} )}E(|d_{|\psi^{-1}_x(y)}\si|^2) dy.$$
Let $z\in B(x,\frac{R+s}{\sqrt d})$ and for every $j\in\{1, \cdots, n\},$ $\si_j\in \rhxd$
be the normed section orthogonal to the hyperplane of sections $\si$ such that $\frac{\partial \si}{\partial y_j|z}=0$.
Still assuming that $dx=dV_h$, we know 
from Lemma 2.1 of \cite{Tian}, see also  Lemma 2.2.3 of \cite{GaWe2}, that
$$ \| \frac{\partial \si_j}{\partial y_j}(z)\|_{h^d}^2 \equid \pi \delta_L d^{n+1},$$
so that again
$$\big| \frac{\partial \si_j}{\partial y_j}(z)\big|^2 \le  \pi \delta_L e^{\pi(R+s)^2} d^{n+1}+o(d^{n+1}),$$
where the $o(d^{n+1})$ can be chosen not to depend on $x\in \R X$.
We deduce that
\beq
E(|d\si_{|z}|^2)&=& \sum_{j=1}^n  E(\big | \frac{\partial \si}{\partial y_j}(z)\big|^2) \\
& = & \sum_{j=1}^n \frac{1}{\sqrt \pi }\int_\R a^2 | \frac{\partial \si_j}{\partial y_j}(z)\big|^2 e^{-a^2} da\\
& \leq &  \frac{n}{2} \pi \delta_L e^{\pi(R+s)^2} d^{n+1}+o(d^{n+1}).
\eeq
Finally,
$  \limsup_{d\to \infty}\sup_{x\in \R X}\frac{1}{d^{n+1}}E\big(\sup_{B(x,\frac{R}{\sqrt d}) }|d\si|^2\big)\leq \frac{n}{2} \pi \delta_L g_R(s).$
By choosing $s$ such that $g_R(s) = \rho_R,$ see (\ref{rhoR}), we obtain the result in the case
where $v= 1$ on $X$. 

In general, the Bergman section  at $x$ for the $L^2$-product (\ref{l2product}) 
associated to the volume form $dx$
is equivalent to the Bergman section $\si_0$ at $x$ for $dV_h$ times $\sqrt {v(x)}$,
because $\si_0$ has its $L^2$-norm concentrated on the ball $B(x,\log d /\sqrt d). $
The same  holds true for 
the $\sigma_j$'s, and  the result follows by replacing $\delta_L$ with $v(x)\delta_L $.
\epr

\subsection{Probability of the local presence of closed affine hypersurfaces}
Following the approach of F. Nazarov and M. Sodin (see \cite{Nazarov-Sodin} or also \cite{Lerario}),
we deduce the following Proposition \ref{prop 5} from Propositions \ref{prop 3} and \ref{prop 4}. 
It estimates from below the probability of presence, in a ball of radius inversely proportional to $\sqrt d$,
of a given affine real algebraic hypersurface in the vanishing locus of sections of high tensor powers of an ample real line bundle.

Let $(X,c_X)$ be a smooth real projective manifold of positive dimension $n$ and $(L,h,c_L)$ be
a real holomorphic Hermitian line bundle of positive curvature over $X$. Let $(U,P)$ be
a regular pair given by Definition \ref{def 1} and $\Sigma= P^{-1}(0)\subset U.$ 
Then, for every $x\in \R X$, we set
 $B_{d}= B(x,\frac{R_\up}{\sqrt d})$ and  denote by 
$Prob_{x,\Sigma}(L^d)$ 
the probability that $\si \in \rhxd$ 
has the property that  $\si^{-1}(0)\cap B_d$
contains a hypersurface $\Sigma'$
 such that the pair $ (B_d,\Sigma') $ be diffeomorphic
 to $(\R^n,\Sigma)$. That is, 
 \beq
  Prob_{x,\Sigma}(L^d) = \mu_\R \big(\{\si \in \rhxd \, | \,
 \si^{-1}(0)\cap B_d \supset \Sigma'  
  \text{ and } (B_d,\Sigma') \sim (\R^n,\Sigma)  \}\big).
 \eeq
We then set $Prob_{\Sigma}(L^d) = \inf_{x\in \R X} Prob_{x,\Sigma}(L^d)$.

\begin{Proposition}\label{prop 5}
Let $X$ be a smooth real projective manifold of positive dimension $n$ and $L$ be
a real holomorphic Hermitian line bundle of positive curvature over $X$. Let $(U,P)$ be
a regular pair given by Definition \ref{def 1} and $\Sigma= P^{-1}(0)\subset U.$ 
Then, 
$$\liminf_{d\to \infty} Prob_{\Sigma}(L^d) \geq m_{\tau_\up}, $$  see (\ref{mtau}).
\end{Proposition}
\bpr
Let $x\in \R X$ and  let us choose a $h$-trivialization of $(L,h)$ given by Definition \ref{trivi}. 
By Proposition \ref{prop 3}, there exist, for every $d$ large enough, a compact $K_d$,
an open set $U_d$ and a section $\si_\up\in \rhxd$ such that 
$$ K_d \subset U_d \subset B_{d}\subset W_x,$$
and $(U_d,\si^{-1}(0)\cap U_d)$ be diffeomorphic to $(U,\Sigma)$.
 Moreover, for every $(\delta,\epsilon)\in \mathcal T_\up$,
\begin{equation}\label{infU1}
\inf_{U_d\setminus K_d} |\si_\up| >\frac{\delta}{2} \sqrt d^n
\end{equation}
and for every $z$ in $U_d$,
\begin{equation}\label{infU2}
  |\si_\up(z)| < \frac{\delta}{2}\sqrt d^n \Rightarrow |d_{|z}\si_\up| > \frac{\epsilon}{2}\sqrt d^{n+1}.
\end{equation}
The moduli $|\si_\up|$ and $|d\si_\up|$ are computed here in the $h$-trivialization of $L^d$, see (\ref{derivee}).	
Denote by $\si_P^\perp$ the orthogonal hyperplane to $\si_\up$
in $\rhxd$ and by $s_P$ the orthogonal symmetry of $\rhxd$
which fixes $\si_P^\perp$.
Then, the average value of $\sup_{B_d } |\theta|^2$ 
on $\si_P^\perp $ satisfies
\beq
E\big(\sup_{B_d } |\theta|^2\big) & = & \int_{\si_p^\perp } \sup_{B_d } |\theta|^2 d\mu_\R(\theta)\\
& = & \int_{\rhxd} \sup_{B_d}\big|\frac{\si + s_P(\si)}{2}\big|^2 d\mu_\R(\si)\\
& \leq & \frac{1}{2}\int_{\rhxd} \sup_{B_d }|\si |^2 d\mu_\R(\si) + 
 \frac{1}{2}\int_{\rhxd} \sup_{B_d }|s_P(\si) |^2 d\mu_\R(\si) \\
 & \leq & E\big(\sup_{B_d }|\si|^2\big)\\
 & \leq & \frac{1}{2} \delta_L  \rho_{R_\up} v(x) d^n+ o(d^n)
\eeq
from Proposition \ref{prop 4}, where the $o(d^{n})$ does not depend on $x\in \R X$.
Likewise, from Proposition \ref{prop 4}, 
\beq
 E\big(\sup_{B_d }|d\theta|^2\big) 
 & \leq & \frac{1}{2} \delta_L \pi n \rho_{R_\up} v(x)d^{n+1}  + o(d^{n+1}),
 \eeq
 where the $o(d^{n+1})$ does not depend on $x\in \R X$. 
From Markov's inequality follows that  for every $M>0$,
$$ \mu_\R \big(\{ \theta \in \si_P^\perp \, | \, \sup_{B_d } |\theta|^2 \geq M^2 v(x) \frac{\delta^2 \delta_L}{4\|P\|_{L^2}^2} d^n\}\big)
\leq \frac{2\|P\|_{L^2}^2 \rho_{R_\up}}{M^2 \delta^2}+ o(1) $$
and 
$$ \mu_\R\big(\{ \theta \in \si_P^\perp \, | \, \sup_{B_d } |d\theta|^2 \geq M^2 v(x)\frac{\epsilon^2 \delta_L}{4\|P\|_{L^2}^2} d^{n+1}\}\big)
\leq \frac{2\pi n \|P\|_{L^2}^2 \rho_{R_\up}}{M^2 \epsilon^2}+ o(1), $$
where $o(1)$ does not depend on $x \in \R X$. 
As a consequence, setting 
\beq
\mathcal E_{ \si_P^\perp} =\big\lbrace \theta \in \si_P^\perp \, | \, 
\sup_{B_d } |\theta|^2 < M^2 v(x)\frac{\delta^2 \delta_L}{4\|P\|_{L^2}^2} d^n 
\text{ and } \sup_{B_d } |d\theta|^2 < M^2 v(x)\frac{\epsilon^2 \delta_L}{4\|P\|_{L^2}^2} d^{n+1}\big\rbrace
\eeq
we have
 \beq
  \mu_\R(\mathcal E_{\si_P^\perp}) \geq 1  - 2\frac{\|P\|_{L^2}^2 \rho_{R_\up}}{M^2 }\big(\frac{1}{\delta^2}+ \frac{\pi n }{\epsilon^2}\big)-  o(1),
\eeq
where the $o(1)$ does not depend on $x$.
Choosing $(\delta,\epsilon) $ which minimizes  the function $(\delta, \epsilon)\mapsto \big(\frac{1}{\delta^2}+ \frac{\pi n }{\epsilon^2}\big) $, 
we deduce from (\ref{tauUP}) that 
$
\mu_\R\big(\mathcal E_{\si_P^\perp}) \geq 1  - \frac{\tau_\up}{M^2 }-  o(1) .
$
Now, setting
$$\mathcal F_{M,\up} = \big\lbrace a \frac{\si_\up}{\|\si_\up\|_{L^2}} + \theta  \in \rhxd \, | \,  
a>M \text{ and } \theta \in \mathcal E_{\si_P^\perp}\big\rbrace,$$
where
\begin{equation}\label{norme}
 \|\si_\up\|_{L^2} \equid \frac{\sqrt{v(x)\delta_L }}{\|p\|} 
\end{equation}
by the first part of Proposition \ref{prop 3} and the fact that the mass of $\si_\up$ concentrates on
small balls $B(x,\log d/\sqrt d)$.  
Take $\si \in \mathcal F_{M,\up}$. 
From the estimates (\ref{infU1}), (\ref{infU2}), (\ref{norme}) and the  definition of $\mathcal E_{\si_P^\perp}$, 
for large enough $d$, $0$ is a regular value of $\si$ and  from Lemma \ref{lemme 5}, $\si^{-1}(0) \cap U_d$
is isotopic to $\si_\up^{-1} (0)\cap U_d$, so that the pair
$\big(B_d,\si^{-1}(0) \cap U_d\big) $ is diffeomorphic to $( B(0,R_\up),\Sigma).$
The result follows from the fact that 
\beq Prob_{x,\Sigma}(L^d)\geq  \mu_\R\big(  \mathcal F_{M,\up} \big) 
&\geq &\big(\frac{1}{\sqrt \pi}\int^{+\infty}_M e^{-t^2}dt \big) \big(1-\frac{\tau_\up}{M^2}-o(1)\big) \\
&= & f_{\tau_\up} (M)-o(1),
\eeq
see \S \ref{para I.1}.
We choose $M\in [\sqrt {\tau_\up}, +\infty[$  which maximizes   $f_{\tau_\up}$, see (\ref{mtau}),
and take the limit.
\epr
\bpr[ of Theorem \ref{theo 1}] Let $\up\in \mathcal I_\Sigma$, see \S \ref{para I.2}. For every $d>0$,
let $\Lambda_d$ be a  maximal  subset of $\R X$ with the property that two distinct points of $\Lambda_d$
are at distance greater than $\frac{2R_\up}{\sqrt d}$. The balls centered at points of $\Lambda_d$
and of radius $\frac{R_\up}{\sqrt d}$ are disjoints, whereas the ones of radius
$\frac{2R_\up}{\sqrt d}$ cover $\R X$. For every $x\in \Lambda_d$ and every $\si \in \rhxd\setminus \R \Delta_d$,
we set $N_\Sigma (x,\si)=1$ if the ball $B_d$ contains a hypersurface
$\Sigma'$ such that $ \Sigma' \subset \si^{-1}(0)$ and $(B_d,  \Sigma')$
is diffeomorphic to $(\R^n, \Sigma)$, whereas $N_\Sigma (x,\si)=0$ otherwise. 
Recall that $N_\Sigma(\si) $ denotes the maximal number of disjoint open subsets 
of $\R X$ having the property that each such open subset $U'$ 
contains a hypersurface $\Sigma'$
such that $\Sigma'\subset \rcsi$
and $(U',\Sigma')$
be diffeomorphic to $(\R^n , \Sigma)$. 
Thus,
\beq
E(N_\Sigma )& \geq & \int_{\rhxd \setminus \R \Delta_d} \big(\sum_{x\in \Lambda_d} N_\Sigma (x,\si) \big) d\mu_\R(\si)\\
& = & \sum_{x\in \Lambda_d} Prob_{x,\Sigma}(L^d)\\
& \geq &|\Lambda_d| Prob_{\Sigma}(L^d)
\eeq
by Proposition \ref{prop 5}. 
We deduce from the inclusion $\R X\subset \bigcup_{x\in \Lambda_d} B\big(x, \frac{2R_\up}{\sqrt d}\big)$
that
\beq
Vol_h(\R X) & \leq & \sum_{x\in \Lambda_d} Vol\big(B(x,\frac{2R_\up}{\sqrt d})\big)\\
& \leq & 2^n | \Lambda_d| Vol\big(B_d\big)+ o\big(\frac{|\Lambda_d|}{\sqrt d^n}\big).
\eeq
From Proposition \ref{prop 5} follows then that 
$$ \liminf_{d\to \infty }  \frac{1}{\sqrt d^n}E(N_\Sigma) \geq \frac{m_{\tau_\up} Vol_h(\R X)}{2^n Vol\big(B_d\big)}.$$
This lower bound  holds for every pair $\up\in \mathcal I_\Sigma$
and we get the result by taking the supremum, see (\ref{csigma}).
\epr
\bpr[ of Corollary \ref{coro 1}] For every $d>0$, 
\beq
E(b_i) &=& \int_{\rhxd \setminus \R \Delta_d} b_i(\R C_\si,\R) d\mu_\R(\si)\\
& \geq & \int_{\rhxd\setminus \R \Delta_d} \big(\sum_{\Sigma \in \mathcal H_n } N_\Sigma (\si) b_i(\Sigma)\big) d\mu_\R(\si)\\
& \geq & \sum_{\Sigma \in \mathcal H_n } b_i(\Sigma) E(N_\Sigma).
\eeq
Hence, the first lower bound follows from Theorem \ref{theo 1}, while the second one 
follows along the same lines. The last part of Corollary \ref{coro 1} is then a consequence of Proposition \ref{prop 2}.
\epr

\bibliographystyle{amsplain}
\bibliography{Lower}

\noindent
\textsc{Universit\'e de Lyon \\
CNRS UMR 5208 \\
Universit\'e Lyon 1 \\
Institut Camille Jordan} \\
43 blvd. du 11 novembre 1918 \\
F-69622 Villeurbanne cedex\\
France

\noindent
gayet@math.univ-lyon1.fr\\
welschinger@math.univ-lyon1.fr

\end{document}